\documentclass[12pt,oneside,reqno]{amsart}
\usepackage{verbatim}
\usepackage[draft]{todonotes}

\usepackage{amsthm}
\usetikzlibrary{arrows}
\usepackage{blindtext}

\usepackage{amssymb}
\usepackage{mathrsfs}
\usepackage[all]{xy}
\usepackage{tikz}
\usepackage[alphabetic,nobysame,abbrev]{amsrefs}
\usepackage{enumerate}
\RequirePackage{fix-cm}
\usepackage{bm}
\usepackage{mdframed}
\usepackage{cancel}
\usepackage{makecell}
\usepackage{diagbox}

\usepackage{amscd}
\usepackage{amsmath}
\usepackage{bm}
\usepackage{newtxtext} 
\numberwithin{equation}{section}
\usepackage{romanbar}
\usepackage[toc,page]{appendix}
\usepackage{hyperref}
\input xy
\xyoption{all}
\usepackage{color}
\RequirePackage{fix-cm}
\usepackage{bm}
\usepackage{mdframed}
\usepackage{cancel}
\usepackage{makecell}
\usepackage{diagbox}
\usepackage{pgf,tikz}
\usepackage{mathrsfs}
\usetikzlibrary{arrows}
\usetikzlibrary{arrows,decorations.markings}
\usetikzlibrary{matrix}
\usetikzlibrary{graphs}
\usetikzlibrary{backgrounds}

\definecolor{qqqqff}{rgb}{0.,0.,1.}
\definecolor{xdxdff}{rgb}{0.49019607843137253,0.49019607843137253,1.}

\definecolor{qqqqff}{rgb}{0.,0.,1.}

\usepackage{tikz}
\usepackage{graphicx}
\usepackage{hyperref}
\usepackage[utf8]{inputenc}
\usepackage{amsmath,amsfonts,amssymb,euscript,graphicx,color,tikz}

\newtheorem{lemma}{Lemma}[section]

\newtheorem{theorem}[lemma]{Theorem}
\newtheorem{proposition}[lemma]{Proposition}

\theoremstyle{definition}
\newtheorem{remark}[lemma]{Remark}
\newtheorem{notation}[lemma]{Notation}
\newtheorem{definition}[lemma]{Definition}

\DeclareMathOperator{\Mod}{Mod}
\DeclareMathOperator{\modd}{mod}

\DeclareMathOperator{\Hom}{Hom}

\DeclareMathOperator{\Ab}{Ab}

\DeclareMathOperator{\Ker}{Ker}
\DeclareMathOperator{\Coker}{Coker}

\DeclareMathOperator{\Mat}{Mat}

\newtheorem*{theorem 0*}{Theorem}
\newtheorem*{theorem a*}{Theorem A}
\newtheorem*{theorem b*}{Theorem B}

\newcounter{diagram}
\numberwithin{diagram}{section}

\begin{document}
	
	\title{Revising Auslander-Gruson-Jensen duality}
	
	\author{Ramin Ebrahimi}
	\address{School of Mathematical Sciences, Zhejiang Normal University, Jinhua 321004, China}
	\email{rebrahimi@zjnu.edu.cn / ramin.ebrahimi1369@gmail.com}
	
	\author{Rasool Hafezi}
	\address{School of Mathematics and Statistics, Nanjing University of Information Science and Technology, Nanjing, Jiangsu 210044, P.R. China}
	\email{hafezi@nuist.edu.cn}

\author{Jiaqun Wei}
	\address{School of Mathematical Science, Zhejiang Normal University, Jinhua 321004, China}
	\email{weijiaqun5479@zjnu.edu.cn}

	

	\begin{abstract}
		For a ring $A$, there is a well-known duality between definable subcategories of right $A$-modules and definable subcategories of left $A$-modules, which is a consequence of Auslander-Gruson-Jensen duality $\modd\text{-}\big(\modd\text{-}A\big)\rightarrow \modd\text{-}\big(\modd\text{-}A^{op}\big)$.
The existence of this duality arises from the fact that $\modd\text{-}\big(\modd\text{-}A\big)$ is the free abelian category over the pre-additive category $A$ with a single object.

		In this note, first, we give a simple description of the free abelian category.
This description clarifies the Auslander-Gruson-Jensen duality and also the duality between the definable subcategories of the right $A$-modules and those of the left $A$-modules.
	\end{abstract}
	
	\maketitle


	\section{Introduction}
    Let $\mathcal{X}$ be a pre-additive category. The {\it free abelian category} over $\mathcal{X}$ is by definition an abelian category $\mathcal{A}$ together with an additive functor $F:\mathcal{X}\rightarrow\mathcal{A}$ such that any other additive functor $\mathcal{X}\rightarrow\mathcal{B}$ to an abelian category $\mathcal{B}$ factors through $F$ by a unique, up to natural isomorphism, functor $\mathcal{A}\rightarrow\mathcal{B}$. The free abelian category always exists and, by the universal property, it is unique up to equivalence.
    
	For a ring $A$, denote by $\Mod\text{-}A$ (resp. $\modd\text{-}A$) the category of all (resp. finitely presented) right $A$-modules. We can consider $A$ as a pre-additive category with a single object. Then, by \cite{G} the free abelian category over $A$ is given by the functor
	\begin{align*}
	H:A\longrightarrow &\big(\modd\text{-}(\modd\text{-}A)\big)^{op}.\\
	A\longmapsto &\Hom_A(A,-)
	\end{align*}
	
Recall that $\modd\text{-}(\modd\text{-}A)$ is the category of all finitely presented covariant functors from $\modd\text{-}A$ to the category of all abelian groups. Note that, since $\modd\text{-}A$ is an additive category with cokernels, by a result of Freyd \cite{Fre65}, $\modd\text{-}(\modd\text{-}A)$ is an abelian category.

The free abelian category $\modd\text{-}(\modd\text{-}A)$ can be used for studying definable subcategories, which are central in model theory of modules \cites{Zie84,Pre09}.
Recall that a subcategory $\mathscr{D}$ of $\Mod\text{-}A$ is called {\it definable} if it is closed under direct products, direct limits and pure subobjects.
It is well-known that a subcategory $\mathscr{D}\subseteq\Mod\text{-}A$ is definable if and only if there is a collection $\mathscr{F}$ of functors in $\modd\text{-}(\modd\text{-}A)$, such that
	\begin{center}
		$\mathscr{D}=\{M\in \Mod\text{-}A\mid \overrightarrow{F}(M)=0, \forall F\in \mathscr{F}\}$,
	\end{center}
	where $\overrightarrow{F}$ is the unique extension of $F$ to $\Mod\text{-}A$ which commutes with direct limits \cite{Pre09}*{Corollary 10.2.32}. So, definable subcategories of $\Mod\text{-}A$ are related to subcategories of the free abelian category $\modd\text{-}(\modd\text{-}A)$.
To get a bijection, Krause proved in \cite{Kra01} that any definable subcategory of $\Mod\text{-}A$ corresponds to a unique Serre subcategory of $\modd\text{-}(\modd\text{-}A)$.

Knowing that $\modd\text{-}\big(\modd\text{-}A\big)$ is the free abelian category over the pre-additive category $A$ with a single object, it is obvious that there exists a canonical exact duality
$\modd\text{-}\big(\modd\text{-}A\big)\rightarrow \modd\text{-}\big(\modd\text{-}A^{op}\big)$, known as Auslander-Gruson-Jensen duality \cites{Au68,GJ}. In the light of this duality and the above mentioned result of Krause, we can see a natural duality between definable subcategories of right $A$-modules and definable subcategories of left $A$-modules.
This result is due to Herzog \cite{Her93} (see also \cite{Baz08}).
	
There is another construction for the free abelian category due to Adelman \cite{Ade73}.
	Adelman's method works for additive categories, not pre-additive categories. So, first we need to consider the category $\mathbb{A}$, the smallest additive category containing the pre-additive category $A$ with a single object. $\mathbb{A}$ is the full subcategory of $\modd\text{-}A$ consisting of all objects of the form $A^n$ for all non-negative integers. For two objects $A^n$ and $A^m$ in $\mathbb{A}$, $\Hom_{\mathbb{A}}(A^n,A^m)=\Mat_{m\times n}(A)$, the set of all $m\times n$ matrices with entries from $A$, and if we consider $U\in \Mat_{m\times n}(A)$ as a morphism $A^n\rightarrow A^m$, it acts from left to column vectors.

$\mathbb{A}$ is an additive category, and we prove that, Adelman's free abelian category over $\mathbb{A}$, denoted by $\Romanbar{3}(\mathbb{A})$, gives us the free abelian category over $A$. Also, we prove that there is an equivalence
	\begin{equation*}
    \big(\modd\text{-}(\modd\text{-}A)\big)^{op}\overset{\sim}{\rightarrow} \Romanbar{3}(\mathbb{A}),
    \end{equation*}
see Theorem \ref{2.7}. Objects of the category $\Romanbar{3}(\mathbb{A})$ are chains of the form $A^n\overset{U}{\rightarrow}A^m\overset{V}{\rightarrow}A^p$, and morphisms are commutative diagrams modulo some homotopy relation, see the next section for more details.
	
	Using this equivalence, we prove the following characterization of definable subcategories of $\Mod\text{-}A$. A subcategory $\mathscr{D}\subseteq\Mod\text{-}A$ is definable if and only if there is a collection of pairs of matrices $(U_i,V_i)_{i\in I}$ (of appropriate size), with entries from $A$, such that 
	\begin{align*}
		\mathscr{D}=\{M\in \Mod\text{-}A \mid \forall i\in I \;\text{and}\; \forall x\in M^m, xU_i=0\implies x=yV_i\; \text{for some}\; y\in M^p\}.
	\end{align*} 
Note that in this equation, $m$ is the number of rows of $U_i$ and $p$ is the number of rows of $V_i$.
This was proved in \cite{E25} as a consequence of a general result about definable subcategories of functor categories.

The presentation of the free abelian category as $\Romanbar{3}(\mathbb{A})$ is much simpler than $\modd\text{-}\big(\modd\text{-}A\big)$ and can be expressed by just matrices of elements of $A$.
For example, the Auslander-Gruson-Jensen duality can be described as follow:
If we denote the Adelman's free abelian category over the opposite ring $A^{op}$ by $\Romanbar{3}(\mathbb{A}^{op})$, there exists a duality $\Romanbar{3}(\mathbb{A})\rightarrow \Romanbar{3}(\mathbb{A}^{op})$, by just reversing all arrows and replacing left multiplication of matrices in column vectors, by right multiplication of matrices in row vectors, see Theorem \ref{2.8}. 
	Also, using the category $\Romanbar{3}(\mathbb{A})$, we can prove the duality between definable subcategories of right $A$-modules and definable subcategories of left $A$-modules due to Herzog \cites{Her93, Baz08}, in an easy and nice way, see Theorem \ref{2.17}.
	
	\section{The results}
Throughout this article, $A$ is a ring with identity. The category of all right $A$-modules is denoted by $\Mod\text{-}A$, and the subcategory of all finitely presented right $A$-modules is denoted by $\modd\text{-}A$. For left $A$-modules we use the notations $\Mod\text{-}A^{op}$ and $\modd\text{-}A^{op}$, where $A^{op}$ stands for the opposite ring of $A$.

	We are interested in definable subcategories of $\Mod\text{-}A$. There are some different characterizations for definable subcategories. Let start from some definitions. As usual, the category of all abelian groups is denoted by $\Ab$.
	\begin{definition}
\begin{itemize}
\item[(1)]
A functor $F:\modd\text{-}A\rightarrow \Ab$ is called {\it finitely presented} if there is a sequence of natural transformations
\[\Hom_A(Y,-)\rightarrow \Hom_A(X,-)\rightarrow F\rightarrow 0\]
with $X,Y\in \modd\text{-}A$, which is exact when evaluated at any object in $\modd\text{-}A$.
\item[(2)]
		An additive functor $F:\Mod\text{-}A\rightarrow \Ab$ is called {\it coherent} if it satisfies one of the following two equivalent conditions.
		\begin{itemize}
			\item[(a)]
			$F$ commutes with direct limits, and it is finitely presented when restricted to $\modd\text{-}A$.
			\item[(b)]
			There is a sequence of natural transformations
			\begin{center}
				$\Hom_A(Y,-)\rightarrow \Hom_A(X,-)\rightarrow F\rightarrow 0$,
			\end{center}
			with $X,Y\in \modd\text{-}A$, which is exact when evaluated at any object in $\Mod\text{-}A$.
		\end{itemize}
\end{itemize}
	\end{definition}
	
	\begin{remark}\label{2.2}
		Any functor $F:\modd\text{-}A\rightarrow \Ab$ can uniquely extend to a functor $\overrightarrow{F}:\Mod\text{-}A \rightarrow \Ab$ which commutes with direct limits. So we can say that coherent functors on $\Mod\text{-}A$ are exactly the unique extensions of finitely presented functors from $\modd\text{-}A$ to $\Ab$, i.e. the unique extensions of objects of $\modd\text{-}(\modd\text{-}A)$.
	\end{remark}
	
	Now we can define definable subcategories. 
	
	\begin{definition}
		A subcategory $\mathscr{D}\subseteq \Mod\text{-}A$ is called {\it definable} if it satisfies one of the following equivalent conditions.
		\begin{itemize}
			\item[(1)] $\mathscr{D}$ is closed under direct products, direct limits and pure subobjects.
			\item[(2)] There is a collection $\mathscr{F}$ of objects of $\modd\text{-}(\modd\text{-}A)$ such that 
			\begin{align*}
				\mathscr{D}=\bigcap_{F\in \mathscr{F}}\Ker(\overrightarrow{F})=\{M\in \Mod\text{-}A \mid \overrightarrow{F}(M)=0, \; \forall F\in \mathscr{F}\}.
			\end{align*}
		\end{itemize}
	\end{definition}
	For the equivalence of these conditions the reader is referred to \cite{Pre09}*{Corollary 10.2.32.}.
	
	\begin{remark}\label{2.4}
		Let $F:\modd\text{-}A\rightarrow \Ab$ be a finitely presented functor with projective presentation
		\begin{equation*}
			\Hom_A(Y,-)\rightarrow \Hom_A(X,-)\rightarrow F\rightarrow 0.
		\end{equation*}
		Then by Yoneda lemma, there is a morphism $f:X\rightarrow Y$ that induces this projective presentation. It is easy to see that a module $M$ belongs to $\Ker(\overrightarrow{F})$ if and only if $M$ is injective with respect to $f$, i.e. for any morphism $g:X\rightarrow M$, there exists a morphism $\bar{g}$ making the following diagram commutative.
		\begin{center}
			\begin{tikzpicture}
				\node (X1) at (0,0) {$X$};
				\node (X2) at (3,0) {$Y$};
				\node (X3) at (0,-2) {$M$};
				\draw [->,thick] (X1) -- (X2) node [midway,above] {$f$};
				\draw [->,thick] (X1) -- (X3) node [midway,left] {$g$};
				\draw [->,thick,dashed] (X2) -- (X3) node [midway,left] {$\bar{g}$};
			\end{tikzpicture}
		\end{center}
	Thus a subcategory $\mathscr{D}$ of $\Mod\text{-}A$ is definable if and only if there is a collection $(f_i:X_i\rightarrow Y_i)_{i\in I}$ of morphisms in $\modd\text{-}A$ such that $\mathscr{D}$ is the subcategory of $\Mod\text{-}A$ consisting of all modules that are injective with respect to all $f_i$'s.	
	\end{remark}

	It is not hard to see that $\modd\text{-}A$ is an additive category with cokernels. Thus, by a result of Freyd, $\modd\text{-}(\modd\text{-}A)$ is an abelian category \cite{Fre65}. Also the canonical embedding $A\rightarrow \big(\modd\text{-}(\modd\text{-}A)\big)^{op}$ given by $A\mapsto\Hom_A(A,-)$, provides the free abelian category over the pre-additive category $A$ with a single object \cite{G}.
	
	In the above, we have defined definable subcategories of $\Mod\text{-}A$ using the free abelian category $\modd\text{-}(\modd\text{-}A)$.
There is another description of the free abelian category due to Adelman \cite{Ade73}.
In the following we briefly recall his construction. For the proofs and more details the reader is refereed to loc.cit.
	
	Let $\mathcal{X}$ be a skeletally small additive category and $\rm Ch^{1,2,3}(\mathcal{X})$ be the category of three term chains (not necessarily complex) over $\mathcal{X}$. We denote a test object of $\rm Ch^{1,2,3}(\mathcal{X})$ by $X_1\overset{f_1}{\rightarrow} X_2\overset{f_2}{\rightarrow} X_3$ and a morphism by a commutative diagram like
	\begin{equation}
		\begin{tikzpicture}
			\node (X1) at (-2,0) {$X_1$};
			\node (X2) at (0,0) {$X_2$};
			\node (X3) at (2,0) {$X_3$};
			\node (X4) at (-2,-2) {$Y_1$};
			\node (X5) at (0,-2) {$Y_2$};
			\node (X6) at (2,-2) {$Y_3$};
			\draw [->,thick] (X1) -- (X2) node [midway,above] {$f_1$};
			\draw [->,thick] (X2) -- (X3) node [midway,above] {$f_2$};
			\draw [->,thick] (X4) -- (X5) node [midway,above] {$g_1$};
			\draw [->,thick] (X5) -- (X6) node [midway,above] {$g_2$};
			\draw [->,thick] (X1) -- (X4) node [midway,left] {$\alpha_1$};
			\draw [->,thick] (X2) -- (X5) node [midway,left] {$\alpha_2$};
			\draw [->,thick] (X3) -- (X6) node [midway,left] {$\alpha_3$};
		\end{tikzpicture}
	\end{equation}
	or simply by the triple $(\alpha_1,\alpha_2,\alpha_3)$. We say that $(\alpha_1,\alpha_2,\alpha_3)$ is null-homotopic if there are morphisms $s:X_2\rightarrow Y_1$ and $t:X_3\rightarrow Y_2$ such that $\alpha_2= g_1s+tf_2$. It is not hard to see that null-homotopic morphisms form an ideal of the additive category $\rm Ch^{1,2,3}(\mathcal{X})$. Following \cite{E25} we denote by $\Romanbar{3}(\mathcal{X})$ the associated factor category. $\Romanbar{3}(\mathcal{X})$ is an abelian category and together with the natural functor $\mathcal{X}\rightarrow \Romanbar{3}(\mathcal{X})$ given by $X\mapsto (0\rightarrow X\rightarrow 0)$ provides the free abelian category over $\mathcal{X}$.
	
	A ring $A$ can be considered as a pre-additive category with a single object, but clearly it is not additive. We want to make it additive, by the minimum of effort.
	Define the additive category $\mathbb{A}$ as follow. $\mathbb{A}$ is the full subcategory of $\modd\text{-}A$ with
	\begin{equation*}
		\rm Obj(\mathbb{A})= \{0, A, A^2, A^3,...\}.
	\end{equation*}
	So $\mathbb{A}$ is an additive category and $\Hom_{\mathbb{A}}(A^n,A^m)=\Mat_{m\times n}(A)$, the set of all $m\times n$ matrices with entries from $A$. Note that we are working with right modules, so if we consider a matrix $U\in \Mat_{m\times n}(A)$ as a right $A$-module homomorphism $A^n\rightarrow A^m$, $U$ acts from left to column vectors in $A^n$.
	
	Let $\mathcal{X}$ be a pre-additive category. By the free additive category over $\mathcal{X}$ we mean an additive category $\mathcal{A}$ together with an additive functor $\mathcal{X}\rightarrow \mathcal{A}$ such that any other additive functor from $\mathcal{C}$ to an additive category $\mathcal{B}$ factor through $\mathcal{C}\rightarrow \mathcal{A}$ by a unique, up to natural isomorphism, additive functor $\mathcal{A}\rightarrow\mathcal{B}$.
	
	\begin{proposition}\label{2.5}
		Considering $A$ as a pre-additive category with a single object, the natural embedding $A\rightarrow \mathbb{A}$, given by $A\mapsto A$, provides the free additive category over $A$.
		\begin{proof}
			Straightforward.
		\end{proof}
	\end{proposition}
	
	We record the following proposition for future references.
	\begin{proposition}\label{2.6}
		\begin{itemize}
			\item[(1)] The natural embedding
			\begin{align*}
				A\longrightarrow &\big(\modd\text{-}(\modd\text{-}A)\big)^{op}\\
				A\longmapsto &\Hom_A(A,-)
			\end{align*}
		     provides the free abelian category over the pre-additive category $A$.
		     \item[(2)] The composition of natural embeddings
		     \begin{align*}
		     	A\longrightarrow \mathbb{A}\longrightarrow&\Romanbar{3}(\mathbb{A})\\
		     	A\longmapsto (0\rightarrow &A\rightarrow 0)
		     \end{align*}
		     provides the free abelian category over the pre-additive category $A$.
		\end{itemize}
	\begin{proof}
		For the proof of $(1)$ see \cite[Lemma 1]{G}. $(2)$ follows from \cite[Theorem 1.14]{Ade73} and Proposition \ref{2.5}. Note that the functor $A\longrightarrow \modd\text{-}(\modd\text{-}A)$ in $(1)$ is contravariant, thus we need the opposite category. 
	\end{proof}
	\end{proposition}
	Using the universal property of the free abelian category, we can prove the following theorem, c.f. \cite[Theorem 2.7]{E25}.
	
	\begin{theorem}\label{2.7}
		There are mutually inverse canonical equivalences
		\begin{center}
			\begin{tikzpicture}
				\node (X1) at (0,0) {$\big(\modd\text{-}(\modd\text{-}A)\big)^{op}$};
				\node (X2) at (4,0) {$\Romanbar{3}(\mathbb{A})$};
				\node (k) at (2,1.3) {$\mathfrak{k}$};
				\node (k') at (2,-0.7) {$\mathfrak{k}^{-1}$};
				\draw [->,thick] (X1) to [out=45,in=135] (X2) node [midway,above] {};
				\draw [->,thick] (X2) to [out=225,in=315] (X1) node [midway,above] {};
			\end{tikzpicture}
		\end{center}
		such that,
		 \begin{center}
			$\mathfrak{k}(F)\cong A^m\oplus A^k\overset{\begin{bmatrix} a&0\\b&-1\end{bmatrix}}{\longrightarrow}A^n\oplus A^k\overset{\begin{bmatrix}c&-d\end{bmatrix}}{\longrightarrow}A^l$,
		\end{center}
		for a fixed projective presentation $\Hom_A(Y,-)\overset{\tilde{f}}{\rightarrow}\Hom_A(X,-)\rightarrow F\rightarrow 0$ for $F$ and fixed free presentations for $X$ and $Y$ as in the following commutative diagram.
		\begin{equation}\tag{Diagram 2.1}
			\begin{tikzpicture}\label{dia2.1}
				\node (X1) at (-2.5,2.5) {$A^m$};
				\node (X2) at (-2.5,0) {$A^k$};
				\node (X3) at (0,2.5) {$A^n$};
				\node (X4) at (0,0) {$A^l$};
				\node (X5) at (2.5,2.5) {$X$};
				\node (X6) at (2.5,0) {$Y$};
				\node (X7) at (5,2.5) {$0$};
				\node (X8) at (5,0) {$0$};
				\draw [->,thick] (X1) -- (X2) node [midway,left] {$b$};
				\draw [->,thick] (X3) -- (X4) node [midway,right] {$c$};
				\draw [->,thick] (X5) -- (X6) node [midway,right] {$f$};
				\draw [->,thick] (X1) -- (X3) node [midway,above] {$a$};
				\draw [->,thick] (X3) -- (X5) node [midway,above] {$p$};
				\draw [->,thick] (X5) -- (X7) node [midway,left] {};
				\draw [->,thick] (X2) -- (X4) node [midway,above] {$d$};
				\draw [->,thick] (X4) -- (X6) node [midway,above] {$q$};
				\draw [->,thick] (X6) -- (X8) node [midway,left] {};
			\end{tikzpicture}
		\end{equation}
	And, $\mathfrak{k}^{-1}(A^m\overset{u}{\rightarrow}A^n\overset{v}{\rightarrow}A^p)=F$, where $F$ admits a projective presentation $\Hom_A(Y,-)\overset{\tilde{f}}{\rightarrow}\Hom_A(X,-)\rightarrow F\rightarrow 0$, where $f:X\rightarrow Y$ fits in a diagram of the following form.
	\begin{equation}\tag{Diagram 2.2}
		\begin{tikzpicture}\label{dia2.2}
			\node (X1) at (-2.5,2.5) {$A^m$};
			\node (X2) at (-2.5,0) {$A^m$};
			\node (X3) at (0,2.5) {$A^n$};
			\node (X4) at (0,0) {$A^p$};
			\node (X5) at (2.5,2.5) {$X$};
			\node (X6) at (2.5,0) {$Y$};
			\node (X7) at (5,2.5) {$0$};
			\node (X8) at (5,0) {$0$};
			\draw [double,-,thick] (X1) -- (X2) node [midway,above] {};
			\draw [->,thick] (X3) -- (X4) node [midway,right] {$v$};
			\draw [->,thick] (X5) -- (X6) node [midway,right] {$f$};
			\draw [->,thick] (X1) -- (X3) node [midway,above] {$u$};
			\draw [->,thick] (X3) -- (X5) node [midway,left] {};
			\draw [->,thick] (X5) -- (X7) node [midway,left] {};
			\draw [->,thick] (X2) -- (X4) node [midway,above] {$vu$};
			\draw [->,thick] (X4) -- (X6) node [midway,left] {};
			\draw [->,thick] (X6) -- (X8) node [midway,left] {};
		\end{tikzpicture}
	\end{equation}
    \begin{proof}
    	For a proof in a more general setting see \cite[Theorem 2.7 and Theorem 2.9]{E25}.
        One can write the proof directly using the following facts:
        \begin{itemize}
        \item[(1)]
        $\big(\modd\text{-}(\modd\text{-}A)\big)^{op}$ and $\Romanbar{3}(\mathbb{A})$ are both the free abelian category over $A$ by Proposition \ref{2.6}.
        \item[(2)]
        By the universal property of the free abelian category, the equivalences $\mathfrak{k}$ and $\mathfrak{k}^{-1}$ exist, and they make the following diagram commutative.
        \begin{center}
			\begin{tikzpicture}
	            \node (X0) at (2,3) {$A$};
				\node (X1) at (0,0) {$\big(\modd\text{-}(\modd\text{-}A)\big)^{op}$};
				\node (X2) at (4,0) {$\Romanbar{3}(\mathbb{A})$};
				\node (k) at (2,1.3) {$\mathfrak{k}$};
				\node (k') at (2,-0.7) {$\mathfrak{k}^{-1}$};
				\draw [->,thick] (X1) to [out=45,in=135] (X2) node [midway,above] {};
				\draw [->,thick] (X2) to [out=225,in=315] (X1) node [midway,above] {};
	             \draw [->,thick] (X0) -- (X1) node [midway,right] {};
	             \draw [->,thick] (X0) -- (X2) node [midway,right] {};
			\end{tikzpicture}
		\end{center}
        \item[(3)]
        By $(2)$ we have that $\mathfrak{k}(\Hom_A(A,-))=0\rightarrow A\rightarrow 0$, and $\mathfrak{k}^{-1}(0\rightarrow A\rightarrow 0)=\Hom_A(A,-)$.
        \item[(4)]
For any object $F\in \big(\modd\text{-}(\modd\text{-}A)\big)^{op}$, if we choose a projective presetation $\Hom_A(Y,-)\overset{\tilde{f}}{\rightarrow}\Hom_A(X,-)\rightarrow F\rightarrow 0$, and free presentations for $X$ and $Y$ as in the \ref{dia2.1}, we are given the following exact commutative diagram in $\big(\modd\text{-}(\modd\text{-}A)\big)^{op}$, where we denote the representable functor $\Hom_A(M,-)$ by just $(M,-)$.
        \begin{center}
        \begin{tikzpicture}
 \node (K2) at (-2,3) {$0$};
 \node (K1) at (0.5,3) {$0$};
 \node (F2) at (-2,1.5) {$(Y,-)$};
 \node (F1) at (0.5,1.5) {$(X,-)$};
 \node (FF) at (3,1.5) {$F$};
 \node (F0) at (4,1.5) {$0$};
 \node (X2) at (-2,0) {$(A^l,-)$};
 \node (X1) at (0.5,0) {$(A^n,-)$};
 \node (Y2) at (-2,-1.5) {$(A^k,-)$};
 \node (Y1) at (0.5,-1.5) {$(A^m,-)$};
 \draw [->,thick] (K2) -- (F2) node [midway,above] {};
 \draw [->,thick] (K1) -- (F1) node [midway,above] {};
 \draw [->,thick] (F2) -- (X2) node [midway,above] {};
 \draw [->,thick] (X2) -- (Y2) node [midway,above] {};
 \draw [->,thick] (F1) -- (X1) node [midway,above] {};
 \draw [->,thick] (X1) -- (Y1) node [midway,above] {};
 \draw [->,thick] (F2) -- (F1) node [midway,above] {};
 \draw [->,thick] (F1) -- (FF) node [midway,above] {};
 \draw [->,thick] (FF) -- (F0) node [midway,left] {};
 \draw [->,thick] (X2) -- (X1) node [midway,left] {};
 \draw [->,thick] (Y2) -- (Y1) node [midway,above] {};
 \end{tikzpicture}
        \end{center}
        Therefore, because we know that $\mathfrak{k}(\Hom(A^r,-)\cong\mathfrak{k}(\Hom(A,-))^r=\big(0\rightarrow A\rightarrow 0\big)^r$ by $(3)$, we can compute $\mathfrak{k}(F)$ using the fact that $\mathfrak{k}$ is an exact functor, and the description of kernel and cokernel in the abelian category $\Romanbar{3}(\mathbb{A})$.
        
        For the construction of $\mathfrak{k}^{-1}$, we use the same argument, and the fact that any object in the Adelman's free abelian category $\Romanbar{3}(\mathbb{A})$, can be constructed from objects of the form $0\rightarrow A^r\rightarrow 0$, and taking kernels and cokernels \cite[Proposition 1.5]{Ade73}.
        \end{itemize}
    \end{proof}
	\end{theorem}

There is a duality $d:\modd\text{-}(\modd\text{-}A)\rightarrow \modd\text{-}(\modd\text{-}A^{op})$ known as Auslander-Gruson-Jensen duality \cites{Au68,GJ}. Here we briefly recall how it acts on objects and morphisms. Let $F$ be an object of $\modd\text{-}(\modd\text{-}A)$ with the projective presentation $\Hom(Y,-)\overset{\Hom(f,-)}{\longrightarrow}\Hom(X,-)\rightarrow F\rightarrow 0$. Then
\[dF=\Ker\big(X\otimes-\overset{f\otimes-}{\longrightarrow} Y\otimes-\big),\]
where karnel is taken in the abelian category $\modd\text{-}(\modd\text{-}A^{op})$. The action of $d$ on morphisms is then obtained by the universal property of kernel. If we consider $A$ and $A^{op}$ as pre-additive categories with a single object, and denote by $(-)^{op}:A\rightarrow A^{op}$ the duality between these two categories, Auslander-Gruson-Jensen duality can be seen as the unique, up to natural isomorphism, exact functor $d$ which makes the following diagram commutative.
           \begin{center}
    			\begin{tikzpicture}
    				\node (X1) at (-3,0) {$A$};
    				\node (X2) at (1,0) {$\modd\text{-}(\modd\text{-}A)$};
    				\node (X3) at (-3,-3) {$A^{op}$};
    				\node (X4) at (1,-3) {$\modd\text{-}(\modd\text{-}A^{op})$};
    				\draw [->,thick] (X1) -- (X2) node [midway,above] {};
    				\draw [->,thick] (X3) -- (X4) node [midway,above] {};
    				\draw [->,thick] (X1) -- (X3) node [midway,left] {$(-)^{op}$};
    				\draw [->,thick] (X2) -- (X4) node [midway,right] {$d$};
    			\end{tikzpicture}
    		\end{center}
In the following theorem we show that, if we use the description $\Romanbar{3}(\mathbb{A})$ for the free abelian category, the Auslander-Gruson-Jensen duality has a simple and beautiful description. Here we denote by $\Romanbar{3}(\mathbb{A}^{op})$ the Adelman's free abelian category for the opposite ring $A^{op}$.

    \begin{theorem}\label{2.8}
    	There exists a duality $d_{\Romanbar{3}}:\Romanbar{3}(\mathbb{A})\rightarrow \Romanbar{3}(\mathbb{A}^{op})$. This duality just reverses all arrows and replaces left multiplication of matrices with column vectors, by right multiplication of matrices with row vectors.
    	\begin{proof}
    		By sending a morphism
    		\begin{center}
    			\begin{tikzpicture}
    				\node (X1) at (-3,0) {$A^n$};
    				\node (X2) at (0,0) {$A^m$};
    				\node (X3) at (3,0) {$A^p$};
    				\node (X4) at (-3,-3) {$A^r$};
    	            \node (X5) at (0,-3) {$A^s$};
    				\node (X6) at (3,-3) {$A^t$};
    				\draw [->,thick] (X1) -- (X2) node [midway,above] {$U.-$};
    				\draw [->,thick] (X2) -- (X3) node [midway,above] {$V.-$};
    				\draw [->,thick] (X4) -- (X5) node [midway,above] {$U'.-$};
    				\draw [->,thick] (X5) -- (X6) node [midway,above] {$V'.-$};
    	            \draw [->,thick] (X1) -- (X4) node [midway,left] {$M.-$};
    				\draw [->,thick] (X2) -- (X5) node [midway,left] {$N.-$};
    				\draw [->,thick] (X3) -- (X6) node [midway,left] {$P.-$};
    			\end{tikzpicture}
    		\end{center}
    	   in $\Romanbar{3}(\mathbb{A})$ to the morphism
    	   \begin{center}
    			\begin{tikzpicture}
    				\node (X1) at (-3,0) {$A^n$};
    				\node (X2) at (0,0) {$A^m$};
    				\node (X3) at (3,0) {$A^p$};
    				\node (X4) at (-3,-3) {$A^r$};
    	            \node (X5) at (0,-3) {$A^s$};
    				\node (X6) at (3,-3) {$A^t$};
    				\draw [->,thick] (X2) -- (X1) node [midway,above] {$-.U$};
    				\draw [->,thick] (X3) -- (X2) node [midway,above] {$-.V$};
    				\draw [->,thick] (X5) -- (X4) node [midway,above] {$-.U'$};
    				\draw [->,thick] (X6) -- (X5) node [midway,above] {$-.V'$};
    	            \draw [->,thick] (X4) -- (X1) node [midway,left] {$-.M$};
    				\draw [->,thick] (X5) -- (X2) node [midway,left] {$-.N$};
    				\draw [->,thick] (X6) -- (X3) node [midway,left] {$-.P$};
    			\end{tikzpicture}
    		\end{center}
    	in $\Romanbar{3}(\mathbb{A}^{op})$ clearly we get a functor, and it is indeed a duality.
    	\end{proof}
    \end{theorem}
    The following proposition is saying that $d_{\Romanbar{3}}:\Romanbar{3}(\mathbb{A})\rightarrow \Romanbar{3}(\mathbb{A}^{op})$ is another description of Auslander-Gruson-Jenson duality.
    
    \begin{proposition}\label{2.9}
    The following diagram of functors commutes, up to natural isomorphism.
     \begin{center}
    			\begin{tikzpicture}
    				\node (X1) at (-4,0) {$\big(\modd\text{-}(\modd\text{-}A)\big)^{op}$};
    				\node (X2) at (0,0) {$\Romanbar{3}(\mathbb{A})$};
    				\node (X3) at (-4,-4) {$\big(\modd\text{-}(\modd\text{-}A^{op})\big)^{op}$};
    				\node (X4) at (0,-4) {$\Romanbar{3}(\mathbb{A}^{op})$};
    				\draw [->,thick] (X1) -- (X2) node [midway,above] {$\mathfrak{k}$};
    				\draw [->,thick] (X3) -- (X4) node [midway,above] {$\mathfrak{k}$};
    				\draw [->,thick] (X1) -- (X3) node [midway,left] {$d$};
    				\draw [->,thick] (X2) -- (X4) node [midway,right] {$d_{\Romanbar{3}}$};
    			\end{tikzpicture}
    		\end{center}
    	\begin{proof}
    	If we consider $A$ as a pre-additive category with a single object, we have functors from $A$ to all of the four categories in the diagram such that all triangle in the following diagram are commutative.
    	\begin{center}
    			\begin{tikzpicture}
    	            \node (X0) at (-2,-2) {$A$};
    				\node (X1) at (-4,0) {$\big(\modd\text{-}(\modd\text{-}A)\big)^{op}$};
    				\node (X2) at (0,0) {$\Romanbar{3}(\mathbb{A})$};
    				\node (X3) at (-4,-4) {$\big(\modd\text{-}(\modd\text{-}A^{op})\big)^{op}$};
    				\node (X4) at (0,-4) {$\Romanbar{3}(\mathbb{A}^{op})$};
    				\draw [->,thick] (X1) -- (X2) node [midway,above] {$\mathfrak{k}$};
    				\draw [->,thick] (X3) -- (X4) node [midway,above] {$\mathfrak{k}^{op}$};
    				\draw [->,thick] (X1) -- (X3) node [midway,left] {$d$};
    				\draw [->,thick] (X2) -- (X4) node [midway,right] {$d_{\Romanbar{3}}$};
    	            \draw [->,thick] (X0) -- (X1) node [midway,above] {};
    				\draw [->,thick] (X0) -- (X2) node [midway,above] {};
    				\draw [->,thick] (X0) -- (X3) node [midway,left] {};
    				\draw [->,thick] (X0) -- (X4) node [midway,right] {};
    			\end{tikzpicture}
    		\end{center}
    	Thus, by the universal property of the free abelian category, $d_{\Romanbar{3}}$ and $\mathfrak{k}^{op}d\mathfrak{k}^{-1}$ are naturally isomorphic. So, $d_{\Romanbar{3}}\mathfrak{k}$ and $\mathfrak{k}^{op}d$ are also naturally isomorphic.
    	\end{proof}
    \end{proposition}
	
	In \cite{E25} we characterized definable subcategories of $\Mod\text{-}\mathcal{X}$ for a skeletally small additive category $\mathcal{X}$. The results can be restricted to modules over a ring. For the sake of completeness and because in this special setting the proofs are simple and concrete, in what follows we will state and prove the main results for this special setting.
	
	Let $\mathscr{D}$ be a definable subcategory of $\Mod\text{-}A$. Then there is a collection $\mathscr{F}$ of functors in $\modd\text{-}(\modd\text{-}A)$, such that
	\begin{equation*}
		\mathscr{D}=\{M\in \Mod\text{-}A \mid \overrightarrow{F}(M)=0, \forall F\in \mathscr{F}\}.
	\end{equation*}
    Given the equivalence $\mathfrak{k}:\big(\modd\text{-}(\modd\text{-}A)\big)^{op}\rightarrow \Romanbar{3}(\mathbb{A})$, we want to describe $\mathscr{D}$ using objects of $\Romanbar{3}(\mathbb{A})$.
	
	\begin{lemma}\label{2.10}
		Let $f:X\rightarrow Y$ be a morphism in $\modd\text{-}A$ and 
	\begin{equation*}
	F=\Coker(\Hom_A(f,-))\in \modd\text{-}(\modd\text{-}A).
	\end{equation*}
	By taking finite free presentations of $X$ and $Y$, We obtain an exact commutative diagram as \ref{dia2.1} in Theorem \ref{2.7}.
		The following conditions are equivalent for an arbitrary object $M\in \Mod\text{-}A$.
		\begin{itemize}
			\item[(1)]
			$M$ is injective with respect to $f$, or equivalently $\overrightarrow{F}(M)=0$.
			\item[(2)]
			For any morphism $\alpha$ that makes the left-hand side of the following diagram commutative, there exists a morphism $\beta$ making the right-hand triangle commutative.
			\begin{center}
				\begin{tikzpicture}
					\node (X1) at (-4,0) {$A^m\oplus A^k$};
					\node (X2) at (0,0) {$A^n\oplus A^k$};
					\node (X3) at (4,0) {$A^l$};
					\node (X4) at (0,-3) {$M$};
					\draw [->,thick] (X1) -- (X2) node [midway,above] {$\begin{bmatrix}a&0\\b&-1\end{bmatrix}$};
					\draw [->,thick] (X2) -- (X3) node [midway,above] {$\begin{bmatrix}c&-d\end{bmatrix}$};
					\draw [->,thick] (X2) -- (X4) node [midway,left] {$\alpha$};
					\draw [->,thick,dashed] (X3) -- (X4) node [midway,right] {$\beta$};
					\draw [->,thick] (X1) -- (X4) node [midway,left] {$0$};
				\end{tikzpicture}
			\end{center}
		\end{itemize}
		\begin{proof}
			Let $M$ be injective with respect to $f$. If $\alpha=[\alpha_1\quad \alpha_2]$ makes the left-hand triangle commutative, then $\alpha_2=0$ and $\alpha_1a=0$. Then by the universal property of cokernel, there is a morphism $g:X\rightarrow M$ with $gp=\alpha_1$, and by $(1)$ there is a morphism $\bar{g}$ with $\bar{g}f=g$. Then it is easy to see that $\beta=\bar{g}q$ has the desired property stated in $(2)$. The proof of the other direction is similar.
		\end{proof}
	\end{lemma}
	
	\begin{notation}\label{2.11}
		Let $\mathbf{A} =A^n\overset{u}{\rightarrow} A^m\overset{v}{\rightarrow}A^p$ be an object in $\Romanbar{3}(\mathbb{A})$. Define
		\begin{center}
			$\Omega_{\mathbf{A}}:=\{M\in \Mod\text{-}A \mid \forall \alpha:A^m\rightarrow M, \text{if}\; \alpha u=0\; \text{then}\; \alpha=\beta v\; \text{for some}\; \beta:A^p\rightarrow M\}$.
		\end{center}
	\end{notation}
	
	\begin{lemma}\label{2.12}
		Let $\mathbf{A}=A^m\overset{u}{\rightarrow}A^n\overset{v}{\rightarrow}A^p$ be an object in $\Romanbar{3}(\mathbb{A})$.
		Then
		\begin{itemize}
			\item[(1)]
			$\Omega_{\mathbf{A}}$ is the subcategory of all modules $M$ such that for any morphism $\alpha$ that makes the left-hand triangle in the following commutative, there exists a morphism $\beta$ that makes the right-hand diagram commutative.
			\begin{center}
				\begin{tikzpicture}
					\node (X1) at (-3,0) {$A^m$};
					\node (X2) at (0,0) {$A^n$};
					\node (X3) at (3,0) {$A^p$};
					\node (X4) at (0,-2) {$M$};
					\draw [->,thick] (X1) -- (X2) node [midway,above] {$u$};
					\draw [->,thick] (X2) -- (X3) node [midway,above] {$v$};
					\draw [->,thick] (X2) -- (X4) node [midway,left] {$\alpha$};
					\draw [->,thick,dashed] (X3) -- (X4) node [midway,right] {$\beta$};
					\draw [->,thick] (X1) -- (X4) node [midway,left] {$0$};
				\end{tikzpicture}
			\end{center}
			\item[(2)]
			By taking cokernel of $v$ and $vu$, if we construct an exact commutative diagram as \ref{dia2.2} in Theorem \ref{2.7}, then $\Omega_{\mathbf{A}}$ coincides the subcategory of all modules that are injective with respect to $f$.
		\end{itemize}
		\begin{proof}
			$(1)$ is Obvious. We want to prove that $(1)$ and $(2)$ are equivalent. By Lemma \ref{2.9} a module $M$ satisfies the desired property in $(2)$ if and only if $M$ has the desired property in $(1)$ for the following diagram.
			\begin{center}
				\begin{tikzpicture}
					\node (X1) at (-4,0) {$A^m\oplus A^m$};
					\node (X2) at (0,0) {$A^n\oplus A^m$};
					\node (X3) at (4,0) {$A^p$};
					\node (X4) at (0,-3) {$M$};
					\draw [->,thick] (X1) -- (X2) node [midway,above] {$\begin{bmatrix}u&0\\1&-1\end{bmatrix}$};
					\draw [->,thick] (X2) -- (X3) node [midway,above] {$\begin{bmatrix}v&-vu\end{bmatrix}$};
					\draw [->,thick] (X2) -- (X4) node [midway,left] {$\alpha$};
					\draw [->,thick,dashed] (X3) -- (X4) node [midway,right] {$\beta$};
					\draw [->,thick] (X1) -- (X4) node [midway,left] {$0$};
				\end{tikzpicture}
			\end{center}
			And we can easily see that this diagram and the diagram in $(1)$ describe the same class of modules like $M$.
		\end{proof}
	\end{lemma}

    \begin{remark}\label{2.13}
    	\begin{itemize}
    		\item[(1)]
    			Let $F\in \modd\text{-}(\modd\text{-}A)$, by Theorem \ref{2.7} and Lemma \ref{2.9}, for any $M\in \Mod\text{-}A$, $M\in \Ker\big(\overrightarrow{F}\big)$ if and only if $M\in \Omega_{\mathfrak{k}(F)}$.
    		\item[(2)]
    		  Let $\mathbf{A}=(A^m\overset{u}{\rightarrow}A^n\overset{v}{\rightarrow}A^p)\in \Romanbar{3}(\mathbb{A})$, by Theorem \ref{2.7} and Lemma \ref{2.11} for any $M\in \Mod\text{-}A$, $M\in \Omega_{\mathbf{A}}$ if and only if $M\in \Ker\big(\overrightarrow{\mathfrak{k}^{-1}(\mathbf{A})}\big)$.  
    	\end{itemize}
    \end{remark}
	
	Now we can prove the following theorem, which gives a characterization of definable subcategories.
	This was proved in \cite{E25} as a consequence of a more general theorem.
	
	\begin{theorem}\label{2.14}
		Let $\mathscr{D}$ be a subcategory of $\Mod\text{-}A$. The following are equivalent.
		\begin{itemize}
			\item[(1)]
			$\mathscr{D}$ is a definable subcategory.
			\item[(2)] There is a collection of pairs of matrices $(U_i,V_i)_{i\in I}$ (of appropriate size) with entries from $A$ such that 
			\begin{align}
				\mathscr{D}=\{M\in \Mod\text{-}A \mid \forall i\in I \;\text{and}\; \forall x\in M^m, xU_i=0\implies x=yV_i\; \text{for some}\; y\in M^p\}.\label{MT}
			\end{align} 
	Where $m$ is the number of rows of $U_i$ and $p$ is the number of rows of $V_i$.
		\end{itemize}
		\begin{proof}
			If $\mathscr{D}$ is a definable subcategory, then there is a family $(f_i:X_i\rightarrow Y_i)_{i\in I}$ of morphisms in $\modd\text{-}A$ such that $\mathscr{D}$ is the subcategory of all modules which are injective with respect to all $f_i$'s. Then by Lemma \ref{2.10} for each $f_i$ we can find a pair of matrices $(U_i,V_i)$ such that a module $M$ is injective with respect to $f_i$ if and only if for each $x\in M^m$, $xU_i=0$ implies that $x=yV_i$ for some $y\in M^p$. This proves $(1)\implies (2)$. The other direction is similar and uses Lemma \ref{2.12} and is left to the reader. 
		\end{proof}
	\end{theorem}
	
	\begin{definition}\label{2.15}
	Let $\mathscr{D}$ be a definable subcategory of $\Mod\text{-}A$, then by Theorem \ref{2.14} there is a collection of pairs of matrices $(U_i,V_i)_{i\in I}$ such that the equality \eqref{MT} is satisfied. If $(U_i,V_i)_{i\in I}$ is the set of all pairs of matrices with this property, we say that $\mathscr{D}$ {\it is defined by the collection} $(U_i,V_i)_{i\in I}$.
	Not that a definable subcategory of left $A$-modules is defined by the collection $(U_i,V_i)_{i\in I}$, means that
	\begin{align*}
				\mathscr{D}=\{M\in \Mod\text{-}A \mid \forall i\in I \;\text{and}\; \forall x\in M^n, U_ix=0\implies x=V_iy\; \text{for some}\; y\in M^m\}.
			\end{align*} 
	\end{definition}
	
Now we want to reinterpret the duality between definable subcategories of right modules and that of left modules.
	We need the following technical proposition.
	
	\begin{proposition}\label{2.16}
		\begin{itemize}
			\item[(1)]
			Let $\mathscr{D}$ be a definable subcategory of $\Mod\text{-}A$. Then, there is a Serre subcategory $\mathcal{S}$ of $\modd\text{-}(\modd\text{-}A)$ such that
			\begin{equation*}
				\mathscr{D}=\bigcap_{F\in \mathcal{S}}\Ker(\overrightarrow{F}).
			\end{equation*}
			And this assignment gives a bijection between Serre subcategories of $\modd\text{-}(\modd\text{-}A)$ and definable subcategories of $\Mod\text{-}A$.
			\item[(2)]
				Let $\mathscr{D}$ be a definable subcategory of $\Mod\text{-}A$. Then, there is a Serre subcategory $\mathscr{S}$ of $\Romanbar{3}(\mathbf{A})$ such that
			\begin{equation*}
				\mathscr{D}=\bigcap_{\mathbf{A}\in \mathscr{S}}\Omega_{\mathbf{A}}.
			\end{equation*}
			And this assignment gives a bijection between Serre subcategories of $\Romanbar{3}(\mathbf{A})$ and definable subcategories of $\Mod\text{-}A$.
			\item[(3)]
			Each of these bijections reverse the inclusion of mentioned subcategories.
		\end{itemize}
	\begin{proof}
		For the proof of $(1)$ see \cite[Section 2.2]{Kra01}. Then $(2)$ follows from $(1)$, Theorem \ref{2.7} and Remark \ref{2.13}.
		$(3)$ is trivial.
	\end{proof}
	\end{proposition}

So, definable subcategories of $\Mod\text{-}A$ are in correspondence with Serre subcategories of the free abelian category $\modd\text{-}(\modd\text{-}A)$ (or $\Romanbar{3}(\mathbf{A})$). Being a duality, Auslander-Gruson-Jensen duality obviously gives a correspondence between Serre subcategoies. Using this correspondence we get a bijective correspondence between definable subcategories of right $A$-modules and that of left $A$-modules.

	\begin{theorem}\label{2.17}
   Let $A$ be a ring.
        \begin{itemize}
        \item[(1)]
        There is an order preserving bijection from definable subcategories of $\Mod\text{-}A$ to definable subcategories of $\Mod\text{-}A^{op}$.
        \item[(2)]
        This bijection sends the definable subcategory $\mathscr{D}$ of $\Mod\text{-}A$ defined by the collection $(U_i,V_i)_{i\in I}$ of pairs of matrices, to the definable subcategory $\bar{\mathscr{D}}$ of $\Mod\text{-}A^{op}$ defined by the collection $(V_i,U_i)_{i\in I}$.
        \end{itemize}
		\begin{proof}
			For any definable subcategory $\mathscr{D}$ of $\Mod\text{-}A$, by Proposition \ref{2.16} there exists a unique Serre subcategory $\mathscr{S}$ of $\Romanbar{3}(\mathbb{A})$ such that
			\begin{equation*}
				\mathscr{D}=\bigcap_{\mathbf{A}\in \mathscr{S}}\Omega_{\mathbf{A}}.
			\end{equation*}
	Apllying Auslander-Gruson-Jensen duality $d_{\Romanbar{3}}$ to $\mathscr{S}$ we obtain a Serre subcategory of $\Romanbar{3}(\mathbb{A}^{op})$.
	If $\mathscr{D}$ is defined by the collection $(U_i,V_i)_{i\in I}$ (corresponds to isoclasses of objects in $\mathscr{S}$), the definable subcategory of $\Mod\text{-}A^{op}$ associated to the Serre subcategory $d_{\Romanbar{3}}(\mathscr{S})$ is defined by the collection $(V_i,U_i)_{i\in I}$ by Theorem \ref{2.8}.
		\end{proof}
	\end{theorem}

	\section*{Acknowledgments}
    Ramin Ebrahimi is supported by Zhejiang Normal University.
	Rasool Hafezi is supported by the National Natural Science of China (Grant No. 12571042).
	Jiaqun Wei is supported by the National Natural Science Foundation of China (Grant Nos. 12571042, 12271249) and the Natural Science Foundation of Zhejiang Province (Grant No. LZ25A010002).
	
	\section*{Declaration}
	The authors declare that they have no conflicts of interest.
	

\end{document}